# An Extension of Tychonoff's Fixed Point Theorem to Quasi-point Separable Topological Vector Spaces

Jinlu Li

Department of Mathematics

Shawnee State University

Portsmouth, Ohio 45662

USA

jli@shawnee.edu

**Abstract**

In this paper, we introduce the concepts of m-quasiconvex, originally m-quasiconvex, and generalized m-quasiconvex functionals on topological vector spaces. Then we extend the concept of point separable topological vector spaces (by the topological dual spaces) to quasi-point separable topological vector spaces by families of generalized m-quasiconvex functionals. We show that every pseudonorm adjoint topological vector space, which includes locally convex topological vector spaces as special cases, is quasi-point separable. By the Fan-KKM theorem, we prove a fixed point theorem on quasi-point separable topological vector spaces. It is an extension of the fixed point theorem on pseudonorm adjoint topological vector spaces proved in [8]. Therefore, it is a proper extension of Tychonoff's fixed point theorem on locally convex topological vector spaces, which are demonstrated by some examples.

1. Introduction

In 1935, Tychonoff generalized the Schauder's fixed-point theorem from Banach spaces to locally convex topological vector spaces. Since then, many authors have extensively studied the properties of locally convex topological vector spaces. Some counter examples of topological vector spaces, which are not locally convex, were provided (see [1−4], [8−9], [12], [14−16]). In 2020, the present author introduced the concept of pseudonorm adjoint topological vector spaces, on which a fixed point theorem is proved. Meanwhile a concrete example was provided in to show that is a proper extension of the Tychonoff's fixed point theorem (see [8]).

Suppose that a pseudonorm adjoint topological vector space $(X, \tau)$ is equipped with a family of $\tau$-continuous quasi-pseudonorms. In [8], we notice that the totality of the family of $\tau$-continuous quasi-pseudonorms plays a very important role for the existence of fixed points of continuous mappings on nonempty convex and compact subsets of $X$. It provides very useful ideas to extend the pseudonorm adjoint topological vector spaces to more general topological vector spaces.

In contrast, of the concept of point separable topological vector spaces by their topological dual spaces, in this paper, we use the term "quasi-point separable" for totality, which means that the spaces are point separable by families of generalized m-quasiconvex functionals. It is a generalization of the concept of point separable topological vector spaces (see section 3).

By using the concept of families of generalized m-quasiconvex functionals, we introduce quasi-point separable topological vector spaces. Similar to the fixed point problems on pseudonorm adjoint topological vector spaces studied in [8], quasi-point separable topological vector spaces have the following properties to ensure the fixed point property (see sections 3 and 4 for details):

(a) In [8] and some other papers, the families of generalized m-quasiconvex functionals that point separates the considered topological vector spaces are also said to be total;
(b) The family of generalized m-quasiconvex functionals is not required to induce (define) the original topology;
(c) In the definition of quasi-point separable topological vector spaces, we have (i) the m-quasiconvex functionals are continuous with respect to the original topology;
(ii) the generalized m-quasiconvex functionals and the weighted functions are not required being continuous;
(d) Every pseudonorm adjoint topological vector space is a quasi-point separable topological vector space;
(e) Every locally convex topological vector space is a quasi-point separable topological vector space.

From (d) above, quasi-point separable topological vector spaces are extensions of pseudonorm adjoint topological vector spaces. Meanwhile, from section 3, we see that the constructions of quasi-point separable topological vector spaces are similar to that of pseudonorm adjoint topological vector spaces. Even though that the extension seems not very significant, but the conditions for a topological vector space to be quasi-point separable are easier checking. Hence, we believe that it is worth to write the new ideas of quasi-point separable topological vector spaces.

By using the Fan-KKM theorem, some authors (see [5−8]) proved the Tychonoff's fixed point theorem, the Brouwer's fixed-point theorem and the Schauder's fixed-point theorem; Park in [10−13] proved some fixed point theorems; and the present author in [8] proved a fixed point theorem on pseudonorm adjoint topological vector spaces. In this paper, we apply the Fan-KKM theorem to prove a fixed point theorem on quasi-point separable topological vector spaces. It is indeed a proper extension of the Tychonoff's fixed point theorem, which is on locally convex topological vector spaces.

2. **Preliminaries**

In this section, we first recall the concepts about pseudonorm adjoint topological vector spaces defined in [8].

**Definition 2.1 [8]** Let $X$ be a vector space with origin $\theta$. A mapping $p: X \to \mathbb{R}^+$ is called a pseudonorm on $X$ if it satisfies the following conditions:

$W_1$. $p(x) \geq 0$, for all $x \in X$ and $p(\theta) = 0$;
$W_2$. $p(-x) = p(x)$, for all $x \in X$;
$W_3$. For any elements $x_1$, $x_2$ of $X$, and $0 \leq \alpha \leq 1$, one has

$$p(\alpha x_1 + (1 - \alpha)x_2) \leq \alpha p(x_1) + (1 - \alpha)p(x_2).$$

**Definition 2.2 [8].** Let $X$ be a vector space. A mapping $q: X \to \mathbb{R}^+$ is called a quasi-pseudonorm on $X$ if there are a pseudonorm $p$ on $X$ and a strictly increasing continuous function $\varphi: \mathbb{R}^+ \to \mathbb{R}^+$ such that

$W_4$. $q(x) \leq \varphi(p(x))$, for all $x \in X$;
$W_5$. $\varphi(0) = 0$.

Here, $q$ is said to be adjoint with the pseudonorm $p$ and the weighted function $\varphi$.

**Definition 2.3 [8].** Let $(X, \tau)$ be a topological vector space. If $X$ is equipped with a family of $\tau$-continuous quasi-pseudonorms $\{q_\lambda\}_{\lambda \in \Lambda}$ associated with a family of $\tau$-continuous pseudonorms $\{p_\lambda\}_{\lambda \in \Lambda}$ and a family of weighted functions $\{\varphi_\lambda\}_{\lambda \in \Lambda}$, then $(X, \tau)$ is called a pseudonorm adjoint topological vector space.

**Definition 2.4 [8].** A family of quasi-pseudonorms $\{q_\lambda\}_{\lambda \in \Lambda}$ equipped on a topological vector space $(X, \tau)$ is said to be total whenever, for $x \in X$, $q_\lambda(x) = 0$ holds, for every $\lambda \in \Lambda$, then it is necessary to have $x = \theta$. A pseudonorm adjoint topological vector space is said to be total if it is equipped with a total family of quasi-pseudonorms.

**Lemma 2.5 [8].** *Every locally convex topological vector space is a pseudonorm adjoint topological vector space.*

The, in rest of this section, we briefly recall the definitions of point separable topological vector spaces and give the definition of m-quasiconvex functionals, by which, we define quasi-point separable topological vector spaces.

Let $(X, \tau)$ be a topological vector space with origin $\theta$ and let $X^*$ denote the dual space of $X$ (the vector space of linear and continuous functionals on $X$). If, for $x \in X$,

$$h(x) = 0, \text{ for all } h \in X^* \quad \text{implies} \quad x = \theta. \tag{2.1}$$

then, $X$ is said to be point separable (by its dual space $X^*$), or $X$ is called a point separable topological vector space.

(2.1) is equivalently defined as follows: $X$ is point separable (by its dual space $X^*$), if and only if, for any distinct elements $x, y \in X$, there is $h \in X^*$ such that $h(x) \neq h(y)$.

It is well known (see [11], or see Theorem 2.18 in [4]):

**Proposition 2.6**. *Every Hausdorff locally convex topological vector space is point separable.*

## 3. M-quasiconvex functionals and quasi-point separable topological vector spaces

Throughout this section, unless otherwise is stated, let $(X, \tau)$ be a topological vector space.

**Definition 3.1.** Suppose that $u$ is a real valued functional defined on $X$. For any $x_1, x_2 \in X$, and $0 \leq \alpha \leq 1$,

(i) if $u$ satisfies that

$$u(\alpha x_1 + (1 - \alpha)x_2) \leq \mathrm{Max}\{u(x_1), u(x_2)\}, \quad (3.1)$$

then, $u$ is said to be quasiconvex.

(ii) if $u$ satisfies that

$$|u(\alpha x_1 + (1 - \alpha)x_2)| \leq \mathrm{Max}\{|u(x_1)|, |u(x_2)|\}, \quad (3.2)$$

then, $u$ is said to be *m-quasiconvex*.

Here the letter m in front of the word quasiconvex means that it is with respect to magnitude. Similar to the property of quasiconvexity, the m-quasiconvexity of $u$ is equivalent to the following fact, for any nonnegative number $a$, the set

$$\{x \in X : |u(x)| \leq a\}$$

is a convex subset of $X$.

**Definition 3.2.** A m-quasiconvex functional $u$ is said to be *originally m-quasiconvex* if $u(\theta) = 0$.

**Definition 3.3**. If $u$ satisfies that, for any $x_1, x_2 \in X$, and $0 \leq \alpha \leq 1$,

$$|u(\alpha x_1 + (1 - \alpha)x_2)| \leq \alpha |u(x_1)| + (1 - \alpha)|u(x_2)|, \quad (3.3)$$

then, $u$ is said to be absolutely convex.

We provide the following two counter examples to show that, neither the quasiconvexity, nor the m-quasiconvexity includes other one.

**Example 3.4**. Let $u$ be a function on $\mathbb{R}$ defined as

$$u(x) = x^2 - 1, \text{ for } x \in \mathbb{R}.$$

Then, $u$ is convex on $\mathbb{R}$, so is quasiconvex. But $u$ is not m-quasiconvex.

**Example 3.5.** Let $u$ be a function on $\mathbb{R}$ defined as

$$u(x) = \begin{cases} 0, & \text{if } -1 \leq x \leq 1, \\ -1, & \text{otherwise.} \end{cases}$$

Then, $u$ is m-quasiconvex. However, $u$ is not quasiconvex.

**Lemma 3.6.** *We list some useful examples of m-quasiconvex functionals.*

(a) *Every absolutely convex functional is m-quasiconvex;*
(b) *Every convex with nonnegative values, or, concave with nonpositive values functional is absolutely convex, so is m-quasiconvex;*
(c) *Every semi-norm is absolutely convex, so is m-quasiconvex.*

*Proof.* The proof is straightforward and it is omitted here.

**Lemma 3.7.** *Suppose that $u$ is an originally m-quasiconvex functional on $X$. Then, for any $x \in X$, the function $|u(tx)|$ is an increasing function with respect to $t > 0$.*

*Proof.* Take arbitrary $0 < t_1 < t_2$. Then

$$|u(t_1 x)| = \left| u\left(\left(1 - \frac{t_1}{t_2}\right)\theta + \frac{t_1}{t_2} t_2 x\right) \right| \leq \text{Max}\{|u(\theta)|, |u(t_2 x)|\} = |u(t_2 x)|. \qquad \square$$

**Definition 3.8.** A functional $v: X \to \mathbb{R}$ is said to be *generalized m-quasiconvex* on $X$ if there are an originally m-quasiconvex convex functional $u$ and a function $\psi: \mathbb{R}^+ \to \mathbb{R}^+$ such that

D$_1$. $|v(x)| \leq \psi(|u(x)|)$, for all $x \in X$;
D$_2$. $\psi(0) = 0$ and $\psi(t) > 0$, for $t > 0$.

Where, $\psi$ is called the weighted function of $v$ associated with the m-quasiconvex functional $u$.

Notice that the originality of the m-quasiconvex functional $u$ on $X$ satisfies $u(\theta) = 0$. Then by conditions D$_1$ and D$_2$, it implies that a generalized m-quasiconvex functional $v: X \to \mathbb{R}$ satisfies

D$_3$. $v(\theta) = 0$.

**Lemma 3.9.** *Let $(X, \tau)$ be a topological vector space with dual space $X^*$. Then, every $u \in X^*$ is*

(a) *$\tau$-continuous originally m-quasiconvex on $X$;*
(b) *$\tau$-continuous generalized m-quasiconvex on $X$, in which the identity function $I_{\mathbb{R}^+}$ on $\mathbb{R}^+$ is the weighted function of $u$ with itself as the associated $\tau$-continuous originally m-quasiconvex functional.*

*Proof.* The proof is straightforward and it is omitted here. $\qquad \square$

**Lemma 3.10.** *Every $\tau$-continuous pseudonorm $p: X \to \mathbb{R}^+$ is $\tau$-continuous originally m-quasiconvex on $X$.*

**Lemma 3.11.** *Let q be a τ-continuous quasi-pseudonorm on X adjoint with a pseudonorm p and associated with a continuous weighted function φ defined in* [8]. *Then q is τ-continuous generalized m-quasiconvex on X.*

*Proof.* The proof is straightforward and it is omitted here. □

**Definition 3.12**. Suppose that a topological vector space $(X, \tau)$ is equipped with a family of generalized m-quasiconvex functionals $\{v_\lambda\}_{\lambda \in \Lambda}$ on X associated with a family of τ-continuous originally m-quasiconvex functionals $\{u_\lambda\}_{\lambda \in \Lambda}$ and a family of weighted functions $\{\psi_\lambda\}_{\lambda \in \Lambda}$. If, for $x \in X$,

$$v_\lambda(x) = 0, \text{ for every } \lambda \in \Lambda, \text{ implies } x = \theta,$$

then $(X, \tau)$ is said to be quasi-point separable.

Where, the family $\{v_\lambda\}_{\lambda \in \Lambda}$, associated with the families $\{u_\lambda\}_{\lambda \in \Lambda}$ and $\{\psi_\lambda\}_{\lambda \in \Lambda}$, is called a quasi-point separating space for this vector space X, which is simply written as $\{v_\lambda, u_\lambda, \psi_\lambda\}_{\lambda \in \Lambda}$. Meanwhile, we say that the vector space X is quasi-point separated by this family of generalized m-quasiconvex functionals $\{v_\lambda, u_\lambda, \psi_\lambda\}_{\lambda \in \Lambda}$.

**Remarks 3.13**. In Definition 3.12, in order to be more general, the generalized m-quasiconvex functionals $\{v_\lambda\}_{\lambda \in \Lambda}$ are not required being τ-continuous on X; the associated family of originally m-quasiconvex functionals $\{u_\lambda\}_{\lambda \in \Lambda}$ are τ-continuous and the weighted functions $\{\psi_\lambda\}_{\lambda \in \Lambda}$ are positive on $(0, \infty)$ and are not required being continuous.

From Lemmas 2.4 and 2.6, we immediately have

**Lemma 3.14**. *Let $(X, \tau)$ be a topological vector space with dual space $X^*$. If X is point separable (by $X^*$), then X is quasi-point separable with the quasi-point separating space $\{u, u, I_{\mathbb{R}^+}\}_{u \in X^*}$ (it is the point separating space).*

*Proof.* The proof follows from Lemma 3.9 immediately. □

From Lemmas 2.6 and 3.14, we immediately have that

**Lemma 3.15**. *Let $(X, \tau)$ be a Hausdorff locally convex topological vector space with dual space $X^*$. Then X is point separable (by $X^*$), so is quasi-point separable with the quasi-point separating space $\{u, u, I_{\mathbb{R}^+}\}_{u \in X^*}$.*

**Lemma 3.16**. *Suppose that a total pseudonorm adjoint topological vector space $(X, \tau)$ is equipped with a total family of τ-continuous quasi-pseudonorms $\{q_\lambda\}_{\lambda \in \Lambda}$ associated with a family of τ-continuous pseudonorms $\{p_\lambda\}_{\lambda \in \Lambda}$ and a family of weighted functions $\{\varphi_\lambda\}_{\lambda \in \Lambda}$. Then X is quasi-point separable with the quasi-point separating space $\{q_\lambda, p_\lambda, \psi_\lambda\}_{\lambda \in \Lambda}$.*

*Proof.* The proof follows from Lemma 3.11 immediately. □

From other point of view, by Lemma 3.16, we repeat Lemma 3.15 as

**Lemma 3.17**. *Suppose that a locally convex topological vector space $(X, \tau)$ is equipped with a total family of seminorms $\{p_\lambda\}_{\lambda \in \Lambda}$ such that the initial topology $\tau$ on $X$ is induced by $\{p_\lambda\}_{\lambda \in \Lambda}$. Then $X$ is quasi-point separable with the quasi-point separating space $\{p_\lambda, p_\lambda, I_{\mathbb{R}^+}\}_{\lambda \in \Lambda}$.*

*Proof.* Notice that the totality of the family of seminorms $\{p_\lambda\}_{\lambda \in \Lambda}$ equipped on the locally convex topological vector space $(X, \tau)$ implies that $(X, \tau)$ is Hausdorff. Then the proof of this lemma follows from Lemmas 2.5 and 3.15 immediately. □

## 4. A fixed-point theorem in quasi-point separable Housdorrf topological vector spaces

**Theorem 4.1.** *Let $(X, \tau)$ be a quasi-point separable Housdorrf topological vector space with a quasi-point separating space $\{v_\lambda, u_\lambda, \psi_\lambda\}_{\lambda \in \Lambda}$. Let $C$ be a nonempty compact convex subset of $X$. Then every continuous mapping from $C$ to itself has a fixed point.*

*Proof.* As mentioned in Remarks 3.13, the generalized m-quasiconvex functionals $\{v_\lambda\}_{\lambda \in \Lambda}$ are not required being $\tau$-continuous on $X$. Meanwhile the associated family of originally m-quasiconvex functionals $\{u_\lambda\}_{\lambda \in \Lambda}$ are $\tau$-continuous and the weighted functions $\{\psi_\lambda\}_{\lambda \in \Lambda}$ are not required being continuous. Therefore, the proof of this theorem is similar to the proof of Theorem 4.1 in [8] with a small difference.

Let $f: C \to C$ be a continuous mapping. The point separating property of $\{v_\lambda\}_{\lambda \in \Lambda}$ implies that a point $x_0 \in C$ is a fixed point of the mapping $f$ if and only if

$$v_\lambda(x_0 - f(x_0)) = 0, \text{ for every } \lambda \in \Lambda.$$

It follows that the mapping $f$ has a fixed point if and only if

$$\bigcap_{\lambda \in \Lambda}\{x \in C : v_\lambda(x - f(x)) = 0\} \neq \emptyset. \tag{4.1}$$

From conditions $D_1$ and $D_2$ in Definition 3.8, it implies that

$$\bigcap_{\lambda \in \Lambda}\{x \in C : u_\lambda(x - f(x)) = 0\} \subseteq \bigcap_{\lambda \in \Lambda}\{x \in C : v_\lambda(x - f(x)) = 0\}.$$

Therefore, to prove (4.1), it is sufficient to show

$$\bigcap_{\lambda \in \Lambda}\{x \in C : u_\lambda(x - f(x)) = 0\} \neq \emptyset. \tag{4.1}'$$

From the compactness of $C$, to prove $(4.1)'$, we only need to prove that the following family

$$\{\{x \in C : u_\lambda(x - f(x)) = 0\} : \lambda \in \Lambda\}$$

has the finite intersection nonempty property. Let $m$ be an arbitrary positive integer and let $\{\lambda_1, \lambda_2, \ldots, \lambda_m\}$ be an arbitrary finite subset of $\Lambda$. To prove $(4.1)'$, by the finite intersection property, we prove

$$\bigcap_{1 \leq k \leq m}\{x \in C : u_{\lambda_k}(x - f(x)) = 0\} \neq \emptyset. \tag{4.2}$$

To prove (4.2), we first prove that, for any $\delta > 0$,

$$\bigcap_{1 \leq k \leq m}\{x \in C : |u_{\lambda_k}(x - f(x))| < \delta\} \neq \emptyset. \tag{4.3}$$

For $\delta > 0$ as given above, assume, on the contrary, that (4.3) does not hold; that is,

$$\text{Max}_{1 \leq k \leq m} |u_{\lambda_k}(x - f(x))| \geq \delta, \text{ for every } x \in C. \tag{4.9}$$

Based on the mapping $f$, we define a set-valued mapping $F: C \to 2^C \setminus \{\emptyset\}$ as follows:

$$F(x) = \{z \in C: \text{Max}_{1 \leq k \leq m} |u_{\lambda_k}(x - f(z))| \geq \delta\}, \text{ for } x \in C.$$

From the hypothesis (4.9), we see that $x \in F(x)$, and therefore, $F(x) \neq \emptyset$, for every $x \in C$. The $\tau$-continuity of $f: C \to C$ and the $\tau$-continuouity of the m-quasiconvex functionals $\{u_\lambda\}_{\lambda \in \Lambda}$ imply that, for every $x \in C$, $F(x)$ is $\tau$-closed. Next, we show that the mapping $F: C \to 2^C \setminus \{\emptyset\}$ is a KKM mapping.

For any given positive integer $n$, take arbitrary $n$ points $x_1, x_2, \ldots, x_n \in C$. For any positive numbers $t_1, t_2, \ldots, t_n$ satisfying $\sum_{i=1}^n t_i = 1$, let $y = \sum_{i=1}^n t_i x_i$. We show that

$$y \in \bigcup_{1 \leq i \leq n} F(x_i). \tag{4.10}$$

Assume, by the way of contradiction, that (4.10) does not hold. Then we have

$$y \notin F(x_i), \text{ for every } i = 1, 2, \ldots, n.$$

It is

$$\text{Max}_{1 \leq k \leq m} |u_{\lambda_k}(x_i - f(y))| < \delta, \text{ for all } i = 1, 2, \ldots, n. \tag{4.11}$$

From assumption (4.9), the inequality (4.11), and by (3.2) in the definition of originally m-quasiconvex functionals $u_{\lambda_k}$, it follows that

$$\begin{aligned}
\delta &\leq \text{Max}_{1 \leq k \leq m} |u_{\lambda_k}(y - f(y))| && \text{(Assumption 4.9)} \\
&= \text{Max}_{1 \leq k \leq m} |u_{\lambda_k}(\sum_{i=1}^n t_i x_i - f(y))| \\
&= \text{Max}_{1 \leq k \leq m} |u_{\lambda_k}(\sum_{i=1}^n t_i (x_i - f(y)))| \\
&\leq \text{Max}_{1 \leq k \leq m} (\text{Max}_{1 \leq i \leq n} |u_{\lambda_k}(x_i - f(y))|) && \text{(by (3.2))} \\
&= \text{Max}_{1 \leq i \leq n} (\text{Max}_{1 \leq k \leq m} |u_{\lambda_k}(x_i - f(y))|) \\
&< \text{Max}_{1 \leq i \leq n} \delta && \text{(by 4.11)} \\
&= \delta.
\end{aligned}$$

It is a contradiction. It implies that $F: C \to 2^C \setminus \{\emptyset\}$ is a KKM mapping with nonempty $\tau$-closed values. Since $C$ is compact, from Fan-KKM Theorem, we have

$$\bigcap_{x \in C} \{z \in C: \text{Max}_{1 \leq k \leq m} |u_{\lambda_k}(x - f(z))| \geq \delta\} = \bigcap_{x \in C} F(x) \neq \emptyset.$$

Then, there is $z_0 \in C$ satisfying

$$\text{Max}_{1 \leq k \leq m} |u_{\lambda_k}(x - f(z_0))| \geq \delta, \text{ for every } x \in C.$$

In particularly, if we take $x = f(z_0)) \in C$, we get

$$0 = \text{Max}_{1 \leq k \leq m} |u_{\lambda_k}(f(z_0) - f(z_0))| \geq \delta.$$

It is a contradiction. So (4.3) is proved. From the continuity of $f$ and the continuity of every $u_{\lambda_k}$, (4.3) implies that

$$\cap_{1 \leq k \leq m} \{x \in C: |u_{\lambda_k}(x - f(x))| \leq \delta\}$$

is a nonempty $\tau$-closed subset of $C$. Take a strictly decreasing sequence of positive numbers $\{\delta_t\}$ with limit 0; i.e., $\delta_t \downarrow 0$, as $t \to \infty$. Then, for $s > t$, we have

$$\cap_{1 \leq k \leq m} \{x \in C: |u_{\lambda_k}(x - f(x))| \leq \delta_s\} \subseteq \cap_{1 \leq k \leq m} \{x \in C: |u_{\lambda_k}(x - f(x))| \leq \delta_t\}.$$

That is, $\{\cap_{1 \leq k \leq m} \{x \in C: |u_{\lambda_k}(x - f(x))| \leq \delta_t\}\}$ is a decreasing (with respect to inclusions) sequence of nonempty closed subsets of the compact set $C$. It follows that

$$\cap_{1 \leq k \leq m} \{x \in C: |u_{\lambda_k}(x - f(x))| = 0\} = \cap_{t=1}^{\infty} \cap_{1 \leq k \leq m} \{x \in C: |u_{\lambda_k}(x - f(x))| \leq \delta_t\} \neq \emptyset.$$

Hence (4.2) is proved. Since $C$ is compact, by the finite intersection property (4.1)' is proved. Then (4.1) follows immediately. □

**Remarks 4.2.** In Theorem 4.1, the quasi-point separable topological vector space $(X, \tau)$ equipped with a quasi-point separating space $\{v_\lambda, u_\lambda, \psi_\lambda\}_{\lambda \in \Lambda}$ needs to be Housdorrf. It is worth to study whether or not that quasi-point separable property of topological vector spaces implies the Housdorrf property. It lefts to interested authors to consider.

In [12], Park studied fixed point problems on point separable topological vector spaces for a class of functions, such as half-continuous functions. Here, we use Theorem 4.1 and Lemma 3.10 to obtain the following fixed point theorem. Since every Hausdorff locally convex topological vector space is a point separable topological vector space, the following theorem is also an extension of the Tychonoff's fixed point theorem to point separable topological vector spaces.

**Theorem 4.3.** *Let $(X, \tau)$ be a point separable topological vector space. Let $C$ be a nonempty compact convex subset of $X$. Then every continuous mapping from $C$ to itself has a fixed point.*

*Proof.* We see that every point separable topological vector space is Hausdorff. From Lemma 3.14 and Theorem 4.1, this theorem follows immediately. □

By Lemma 3.15 and Theorem 4.1, consequently, we obtain Theorem 4.1 in [8].

**Theorem 4.1 [8].** *Let $(X, \tau)$ be a Housdorrf and total pseudonorm adjoint topological vector space. Let $C$ be a nonempty compact convex subset of $X$. Then every continuous mapping from $C$ to itself has a fixed point.*

From Lemma 3.8 [8] and Theorem 4.1, we obtain

**Tychonoff's fixed point theorem** [12]: *Let X be a Hausdorff locally convex topological vector space. For any nonempty compact convex set C in X, any continuous function f : C → C has a fixed point.* □

## 5. An example: quasi-point separable is a proper extension of locally convexity

In 1935, Tychonoff proved that the topological vector space $l_r$, for $0 < r < 1$, is not locally convex (see [4], [8−9], [17]). In [8], it is proved that the topological vector space $l_r$, for $0 < r < 1$, is a Hausdorff pseudonorm adjoint topological vector space.

In this section, one-step further, we give a simple proof showing that $l_r$, is point separable; and from Lemma 3.9 or 3.14, it follows immediately that $l_r$, for $0 < r < 1$, is quasi-point separable.

Then, from Lemma 3.8 in [8] and Lemma 3.15, every Hausdorff locally convex topological vector space is a Hausdorff quasi-point separable topological vector space. Hence, the fixed point Theorem 4.1 also properly extends Tychonoff fixed point theorem. In this section, Let $S$ denote the set of real sequences.

**Example 5.1**. For every given $r \in (0, 1)$, define a subspace $l_r$ of $S$ as below

$$l_r = \{\{x_i\} \in S : \sum_{i=1}^{\infty} |x_i|^r < \infty\}.$$

Define a functional $q$ on $l_r$ as

$$q(\{x_i\}) = \sum_{i=1}^{\infty} |x_i|^r, \text{ for } \{x_i\} \in l_r.$$

It is known that the functional $q$ induces a metric on $l_r$ by

$$q(\{x_i\} - \{y_i\}) = \sum_{i=1}^{\infty} |x_i - y_i|^r, \text{ for every } \{x_i\}, \{y_i\} \in l_r. \tag{5.1}$$

Let $l_r^*$ denote the topological dual space of $l_r$. For every fixed $n = 1, 2, \ldots$, define a functional $u_n$ on $l_r$ as follows:

$$u_n(\{x_i\}) = x_n, \text{ for } \{x_i\} \in l_r. \tag{5.2}$$

Then, the metric vector space $l_r$ is point separable by $l_r^*$. Furthermore, $l_r$ is quasi-point separable with a quasi-point separating space $\{u_n, u_n, I_{\mathbb{R}^+}\}_{n \in \mathbb{N}}$.

*Proof.* From (5.2), for $n = 1, 2, \ldots$, $u_n$ is a linear and continuous functional on $l_r$. So $\{u_n\} \subseteq l_r^*$. We see that, for $\{x_i\} \in l_r$,

$$u_n(\{x_i\}) = 0, \text{ for all } n = 1, 2, \ldots, \quad \text{implies} \quad \{x_i\} = \theta.$$

Hence, the space $l_r$ is point separated by $\{u_n\}$. Therefore, $l_r$ is point separable (by $l_r^*$).

**Example 5.2**. For every given $r \in (0, 1)$, define a subspace $l^r$ of $S$ as below

$$l^r = \left\{\{x_i\} \in S : \sum_{i=1}^{\infty} \frac{|x_i|^r}{1+|x_i|^r} < \infty\right\}.$$

Define a functional $p$ on $l^r$ as

$$p(\{x_i\}) = \sum_{i=1}^{\infty} \frac{|x_i|^r}{1+|x_i|^r}, \text{ for } \{x_i\} \in l^r.$$

One can show that the functional $p$ induces a metric on the space $l^r$ as follows:

$$p(\{x_i\} - \{y_i\}) = \sum_{i=1}^{\infty} \frac{|x_i - y_i|^r}{1+|x_i - y_i|^r}, \text{ for every } \{x_i\}, \{y_i\} \in l^r. \qquad (5.3)$$

Let $l^{r*}$ denote the topological dual space of $l^r$. Then, the metric vector space $l^p$ is point separable (by $l^{r*}$); so is quasi-point separable with a quasi-point separating space $\{u_n, u_n, I_{\mathbb{R}^+}\}_{n \in \mathbb{N}}$.

*Proof.* To show that the functional defined in (5.3) defines a metric on $l^p$, one may consider the following steps:

(a) $f(t) = \frac{t^r}{1+t^r}$ is a strictly increasing function on $\mathbb{R}^+$;
(b) $(t+s)^r \leq t^r + s^r$, for $t, s \in \mathbb{R}^+$ (this should be proved in Example 5.1);
(c) $\frac{(t+s)^r}{1+(t+s)^r} \leq \frac{t^r}{1+t^r} + \frac{s^r}{1+s^r}$, for $t, s \in \mathbb{R}^+$.

Rest of the proof is almost the same to the proof of Example 5.1. For every $n = 1, 2, \ldots$, define a functional $u_n$ on $l^r$ as in (5.2). The space $l^r$ is point separated by $\{u_n\}$. From $\{u_n\} \subseteq l^{r*}$, $l^p$ is point separated by $l^{r*}$ and $l^r$ is point separable, so is quasi-point separable.

### Acknowledgements

The author sincerely thanks Professor Robert Mendris for his kind communications and comments about this paper.

### References

[1] L.E. J. Brouwer, Uber Abbildung von Mannigfaltigkeiten, Math. Ann. 7l (1912) 97–115.

[2] T. Butsan, S. Dhompongsa, W. Fupinwong, Schuader's conjecture on convex metric spaces, J. Nonlinear Convex Anal. 11, no. 3, (2010) 527–535.

[3] R. Cauty in Solution du problème de point fixe de Schauder, Fund. Math. 170 (2001) 231–246].

[4] K. Conrad, "$L_p$ spaces for $0 < p < 1$", https://kconrad.math.uconn.edu/blurbs/topology/finite-dim-TVS.pdf.

[4] M. M. Day, "The spaces Lp with $0 < p < 1$," Bull. Amer. Math. Soc. 46 (1940), 816–823. Online at https://projecteuclid.org/download/pdf 1/euclid.bams/1183503236.

[5] K. Fan, A generalization of Tychonoff's fixed point theorem, *Math. Ann.*, **142** (1961) 305–310.

[6]     Olga Hadžić, A fixed point theorem in topological vector spaces, Review of Research Faculty of Science, University of Novi Sad, Volume 10 (1980).

[7]     Won Kyu Kim, A fixed point theorem in a Hausdorff topological vector space, Comment. Math. Univ. Carolin. 36, 1 (1995)33–38.

[8]     Jinlu Li, An Extension of Tychonoff's fixed point theorem to pseudonorm adjoint topological vector spaces, *Optimizations*. Published online: 07 Jul 2020 https://doi.org/10.1080/02331934.2020.1789639.

[9]     W. Rudin, *Functional Analysis*. International Series in Pure and Applied Mathematics. 8 (Second ed.). New York, NY: McGraw-Hill Science/Engineering/Math. ISBN 978-0-07-054236-5. OCLC 21163277 (1991).

[10]    Sehie Park, The KKM principle implies many fixed point theorems, Topology and its Applications 135 (2004) 197–206

[11]    S. Park, Recent applications of Fan-KKM theorem, Lecture notes in Math. Anal. Kyoto University, **1841**, 58–68 (2013).

[12]    S. Park, A generalization of the Brouwer fixed point theorem, Bull. Korean Math. Soc. 28 (1991), 33–37.

[13]    S. Park, Ninety years of the Brouwer fixed point theorem, Vietnam J. Math. 27 (1999), 187–222.

[14]    J. Schauder, *Der Fixpunktsatz in Funktionalräumen*, Studia Math. 2 (1930) 171–180

[15]    Joel H. Shapiro, The Schauder Fixed-Point Theorem, An Infinite Dimensional Brouwer Theorem, A Fixed-Point Farrago pp 75–81

[16]    I. Termwuttipong and T. Kaewtem, Fixed point theorems for half-continuous mappings on topological vector spaces, Fixed Point Th. Appl. 2010 (2010), Article ID 814970.

[17]    A. Tychonoff, *Ein Fixpunktsatz*, Mathematische Annalen 111 (1935) 767–776